\documentclass{amsart}
\usepackage{array}
\usepackage{supertabular}

 \newcommand{\Pn}{\mathbb{P}^n}
 
 \renewcommand{\O}{\mathcal{O}}
 
 \DeclareMathOperator{\mult}{mult}
 
 \DeclareMathOperator{\indet}{indet}
 
 \DeclareMathOperator{\inverse}{^{-1}}

 \DeclareMathOperator{\Spec}{Spec}

 \DeclareMathOperator{\Pic}{Pic}
 \DeclareMathOperator{\Hom}{Hom}

 \DeclareMathOperator{\ord}{ord}

 \renewcommand{\P}{\mathbb{P}}
 
 \newcommand{\D}{\mathbb{D}}

\begin{document}

 \newtheorem{thm}{Theorem}[section]
 \newtheorem{cor}[thm]{Corollary}
 \newtheorem{lem}[thm]{Lemma}
 \newtheorem{prop}[thm]{Proposition}

 \theoremstyle{definition}
 \newtheorem{defn}[thm]{Definition}
 \newtheorem{exam}[thm]{Example}

 \theoremstyle{remark}
 \newtheorem{obs}[thm]{Observation}
 \newtheorem{rem}[thm]{Remark}
 \newtheorem{fact}[thm]{Fact}
 \newtheorem{notation}[thm]{Notation}
 \numberwithin{equation}{section}

%--------------------------------------------------

\fontsize{10}{12}\selectfont

\title[]
 {Morphisms from Quintic Threefolds to Cubic Threefolds are Constant}

\author{ David C. Sheppard  }

\address{Department of Mathematics, MIT }

\email{sheppard@math.mit.edu}

\begin{abstract}
We show that every morphism from a quintic threefold in $\mathbb
P^4$ to a nonsingular cubic threefold in $\mathbb P^4$ is constant
in characteristic zero.  In the process, we classify morphisms
from $\mathbb P^2$ to nonsingular cubic hypersurfaces in $\P^4$
given by degree 3 polynomials.
\end{abstract}

%%% ----------------------------------------------------------------------
\maketitle
%%% ----------------------------------------------------------------------
\section*{Introduction}

The author shows in \cite{Sheppard} that if $f:X\rightarrow Y$ is
a morphism of hypersurfaces in $\P^4$ such that $\deg Y=3$ and
$\deg X\leq 4$, then $f$ is either constant or $\deg X=3$ and $f$
is an isomorphism. The purpose of this paper is to extend this
result by proving the following theorem.

\subsection*{Theorem 1}\emph{
If $f:X\rightarrow Y$ is a morphism of hypersurfaces in $\P^4$
over an algebraically closed field of characteristic zero such
that $\deg X=5$, $\deg Y=3$, and $Y$ is nonsingular, then $f$ is
constant.} \vspace{.35cm}

The motivation for investigating such morphisms to cubic
hypersurfaces is the expectation that if $f:X\rightarrow Y$ is a
morphism between hypersurfaces in $\P^n$, $n\geq 4$, such that $Y$
is nonsingular of degree at least 2, then $\deg Y$ divides $\deg
X$ with quotient $q$, and $f$ is given by polynomials of degree
$q$, i.e. $f^*\O_Y(1)=\O_X(q)$.  This result is proven in
\cite{Sheppard} when $\deg Y\geq n+2$.  It is also proven in some
cases where $3\leq \deg Y\leq n+1$. The fact that morphisms from
quintic threefolds to cubic threefolds are necessarily constant is
the first nontrivial case of morphisms to cubics.

\section{Outline of Proof}

Throughout the paper, the base field will be algebraically closed
of characteristic zero, and $f:X\rightarrow Y$ will denote a
nonconstant morphism of hypersurfaces in $\P^4$, such that $\deg
X=5$, $\deg Y=3$, and $Y$ is nonsingular.

The Grothendieck-Lefschetz Theorem, \cite[Theorem 4.3.2]{H1},
states that $\Pic X$ is generated by $\O_X(1)$.  So
$f^*\O_Y(1)=\O_X(m)$ for some positive integer $m$. In terms of
$m$, \cite[Theorem 1]{Sheppard} states that $\deg
f^*c_\text{3}\left(\Omega^1_Y(2)\right)\leq \deg
c_\text{3}\left(\Omega^1_X(2m)\right)$.  By computing the Chern
classes in this inequality, one checks that $m\leq 3$.

We claim there is a rational map $F:\P^4\dashrightarrow\P^4$
undefined at only finitely many points disjoint from $X$ such that
$f=F|_X$.  To see this, consider the short exact sequence
 $$
 H^0(\P^4,\O_{\P^4}(m))\longrightarrow H^0(X,\O_X(m))
 \longrightarrow H^1(\P^4, \O_{\P^4}(m-5))
 $$
Since $H^1(\P^4, \O_{\P^4}(m-5))=0$, the sections of $\O_X(m)$
that define $f:X\rightarrow \P^4$ lift to $\P^4$ and define a
rational map $F:\P^4\dashrightarrow\P^4$, as claimed.

We introduce a hypersurface $H$ in $\P^4$ that will play a key
role.   Lazarsfeld proves in \cite{L} that $\P^k$ is the only
smooth $k$-dimensional variety that is the image of a morphism
from $\P^k$.  In particular, $Y$ is not the image of a morphism
from $\P^3$.  It follows by considering a general hyperplane in
$\P^4$ that the image of $F:\P^4\dashrightarrow\P^4$ is not $Y$.
So $F$ is dominant because the image of $F$ is irreducible and
contains $Y$. Hence $F^{-1}(Y)=X+H$ for some hypersurface $H$ of
degree $3m-5$.

This shows $m\neq 1$.  We see that $m\neq 2$ as follows.  Suppose
$m=2$.  Then $H$ is a hyperplane in $\P^4$.  Therefore
$F|_H:H\dashrightarrow Y$ is not a morphism.  Let $p\in H$ be a
point of indeterminacy of $F$.  This means that $F=(F_0,\dots,
F_4)$ for homogeneous polynomials $F_i$ of degree 2 that all
vanish at $p$. If $Y$ is defined by the homogeneous cubic
polynomial $G=G(y_0,\dots,y_4)$, then $F^{-1}(Y)=X\cup H$ is
defined by the homogeneous sextic polynomial $G(F_0,\dots,F_4)$.
Since the $F_i$ all vanish at $p$ and $G$ is homogeneous of degree
3, $G(F_0,\dots,F_4)$ vanishes to order at least 3 at $p$.  In
other words, $p$ is a triple point of $F^{-1}(Y)=X\cup H$.  This
is impossible because $p$ is contained in the hyperplane $H$ but
not contained in $X$.  Therefore, $m\neq 2$.

It is considerably more difficult to check that $m=3$ is also not
possible.  That task will occupy us for the remainder of the
paper.  We focus our attention on the degree 4 hypersurface $H$
and the rational map $F|_H:H\dashrightarrow Y$.

We prove Theorem 1 by considering the various possibilities for
$H$.  If $H$ contains many copies of $\P^2$, then we will use the
resulting maps $\P^2\rightarrow Y$, which would be given by degree
3 polynomials. So our first task is to classify such morphisms
from $\P^2$.  This is done in Section 2.  Then we return to the
map $H\dashrightarrow Y$ in Section 3 and work through several
cases using the geometry of rational maps from threefolds to three
dimensional cubics.

\section{Morphisms from $\P^2$ to Cubic Threefolds}

Let $g:\P^2\rightarrow Y$ be a morphism given by degree 3
polynomials, i.e. $g^*\O_Y(1)=\O_{\P^2}(3)$.  The purpose of this
section is to prove the following result about $g$.

\subsection*{Theorem 2} \emph{We can choose coordinates\ $x_0,
x_1, x_2$ on $\P^2$ and coordinates on $\P^4$ such that
$$
g=(x_0^3,\ x_1^3 ,\ x_2^3,\ x_0x_1x_2,\ 0):\P^2\longrightarrow
Y\subset\P^4.
$$}
\vspace{.35cm}

Let $S$ be the image of $g$.  Since $Y$ is nonsingular, $S\subset
Y$ is a Cartier divisor.  So by the Grothendieck-Lefschetz
Theorem, \cite[Theorem 4.3.2]{H1}, $S$ is the zero locus $V(s)$ of
a section $s\in H^0(Y,\O_Y(a))$ for some $a>0$. To see that $S$ is
the scheme-theoretic intersection of $Y$ with another hypersurface
$Y^\prime$ in $\P^4$, consider the following piece of a long exact
sequence of cohomology groups.
 $$
  H^0(\P^4, \O_{\P^4}(a))\longrightarrow H^0(Y, \O_Y(a))\longrightarrow
  H^1(\P^4, \O_{\P^4}(a-3))
 $$
Since $H^1(\P^4, \O_{\P^4}(a-3))=0$, $s$ is the image of some
$s^\prime\in H^0(\P^4, \O_{\P^4}(a))$.  So $S=V(s^\prime)\cap Y$.
Take $Y^\prime:=V(s^\prime)$.

Calculate $g^* c_1(\O_Y(1))^2$ to be $9=\deg g\cdot \deg S$ to see
that $3=\deg g\cdot\deg Y^\prime$.  Hence we will break the proof
of Theorem 2 into two cases depending on whether $\deg Y^\prime=3$
or $\deg Y^\prime = 1$.  First, we recall the following result
from \cite{CG}, which we will use repeatedly.

\begin{prop}
Let $\Delta$ be the space of lines on the smooth cubic threefold
$Y$.  Then $\Delta$ is a complete nonsingular surface that does
not contain a rational curve.
\end{prop}

\subsection*{Case 1}

Assume $Y^\prime$ has degree 3.  We will derive a contradiction.

By assumption, $g:\P^2\rightarrow S$ has degree 1, so it is a
finite birational morphism. Our strategy is to analyze the double
point class $\D(g)$ of $g:\P^2\rightarrow Y$.  See \cite[Section
9.3]{F} for the construction and computation of $\D(g)$. The
construction will be used implicitly in the proof of Proposition
2.6 and Corollary 2.7. Following \cite{F}, we calculate the double
point class of $g$:
\begin{align*}
\D(g)  & = g^*g_*[\P^2]-c_1(g^*T_Y)+c_1(T_{\P^2})       \\
       & = g^*c_1(\O_Y(3))-g^*c_1(\O_Y(2))+c_1(\O_{\P^2}(3)) \\
       & = c_1(\O_{\P^2}(6)).
\end{align*}

An important fact for our purposes is that the construction of the
cycle $\D(g)$ not only gives a 1-cycle modulo rational
equivalence, it actually constructs a Weil divisor in $\P^2$. We
will denote this Weil divisor by $D(g)$ and consider it as a
closed subscheme of $\P^2$. This notation differs slightly from
\cite{F} in that we use $D(g)$ to denote a divisor, not just a
set. Roughly speaking, $D(g)$ is the curve in $\P^2$ consisting of
the closed points $x$ such that either $g(x)=g(y)$ for some $y\neq
x$ or such that $T_x\P^2\rightarrow T_{g(x)}Y$ is not injective.
The scheme structure of $D(g)$ comes from the fact that an
integral curve $D$ in $D(g)$ will appear with multiplicity if for
a general point $x\in D$ there is more than one other point $y\in
\P^2$ with $g(x)=g(y)$ or if $g$ ramifies to high order along $D$.

The following result tells us that the image of $D(g)$ under $g$
is equal to the singular locus of $S$ as a set.

\begin{lem}
Let $g:V\rightarrow W$ be a finite surjective birational morphism
of varieties with $V$ regular.  Then for $w\in W$ a closed point,
$w$ is a nonsingular point of $W$ if and only if the
scheme-theoretic preimage $g^{-1}(w)$ is a single reduced point.
\begin{proof}
Suppose $g^{-1}(w)$ is a single reduced point $v$.  Let $D\subset
W$ be a general curve in $W$ through $w$ such that $D$ is a local
complete intersection in $W$ at $w$ and
$$
\dim T_wW-\dim T_wD=\dim W-1.
$$
By assumption, the pullback of the maximal ideal of $w$ generates
the maximal ideal of $v$. So the curve $C:=g^{-1}(D)$ is
nonsingular at $v$ because $C$ is cut out near $v$ by the $\dim
V-1$ equations that define $D$ near $w$. The morphism $g|_C$ is an
isomorphism in a neighborhood of $v$ because $g^{-1}(w)=v$ and
$g|_C$ is birational in a neighborhood of $v$. So $D$ is regular
at $w$, whence $W$ is too by the above equation.

Conversely, suppose $w $ is a regular closed point of $W$.  Let
$A$ be the ring of regular functions on an open affine $U$ of $W$
containing $w$, and let $B$ be the ring of regular functions on
the open affine $g^{-1}(U)$.  Let $m\subset A$ be the maximal
ideal corresponding to the point $w$, and let $M$ denote the
multiplicative set $A\setminus m$. Then $M^{-1} B$ is the integral
closure of $A_m=M^{-1} A$ in its field of fractions. But $A_m$ is
integrally closed because it is regular. Therefore $M^{-1} B=A_m$,
whence $M^{-1} B$ is a local ring whose maximal ideal is generated
by $m$.  In other words, $g^{-1} (w)$ is a single reduced point.
\end{proof}
\end{lem}

Now let $C\subset S$ denote the image of $D(g)$ with its reduced
structure:
$$
C : = g(D(g))_{\text{red}}
$$
By Lemma 2.2, $C$ is the maximal reduced curve contained in the
singular locus of $S$.  We will rule out the case $\deg Y^\prime
=3$ by comparing $D(g)$ to $g^{-1}(C)$, where $g^{-1}(C)$ denotes
the scheme-theoretic preimage.  So our strategy is to compare the
multiplicity of integral curves in $D(g)$ and $g^{-1}(C)$.

\begin{lem}
There is a cubic hypersurface $Y^{\prime\prime}$ in $\P^4$ such
that $S=Y\cap Y^{\prime\prime}$, and $Y^{\prime\prime}$ is
singular along $C$.
\begin{proof}
Let $\{Y_t\}$ be the pencil of cubic hypersurfaces spanned by $Y$
and $Y^\prime$.  Since $Y_t\cap Y=S$ so long as $Y_t\neq Y$, it
suffices to show that one of the $Y_t$ is singular along $C$.

Let $C_0$ be an irreducible component of $C$.  For every point
$p\in C_0$ we have $T_pS=T_pY$ because $S$ is singular at $p$.
Therefore $T_pS=T_pY_t$ for $t$ general, and there is a unique
$t$, say $t_p$, such that $Y_{t_p}$ is singular at $p$.  Only
finitely many of the $Y_t$ are singular, so one of the $Y_t$, say
$Y_{t_0}$, is singular at infinitely many points of $C_0$.  Hence
$Y_{t_0}$ is singular along $C_0$.

If $C_1$ is another irreducible component of $C$, then $C_0$ and
$C_1$ meet at some point $p$ because $C_0$ and $C_1$ are both
images of curves in $\P^2$.  If $Y_{t_1}$ is singular along $C_1$,
then $Y_{t_0}$ and $Y_{t_1}$ are both singular at $p$, so
$Y_{t_0}=Y_{t_1}$.  The Lemma follows.
\end{proof}
\end{lem}

We now assume $Y^\prime$ to be singular along $C$.

\begin{lem}
The scheme-theoretic preimage $g^{-1}(C)$ does not contain a curve
of degree 6 or more.
\begin{proof}
Let $G^\prime=G^\prime(y_0,\dots,y_4)$ be the homogeneous equation
for $Y^\prime$ in $\P^4$.  Set $C_i:=S\cap V(\frac{\partial
G^\prime}{\partial y_i})$ and $D_i:=g^{-1}(C_i)$ for
$i=0,\dots,4$. The $D_i$ are all degree 6 plane curves containing
$g^{-1}(C)$. So it suffices to show that the $D_i$ are not all
equal.

Suppose they were.  Then the equations for the $D_i$ in
$H^0(\P^2,\O_{\P^2}(6))$ are all scalar multiples of each other.
These equations are the images of the partials $\frac{\partial
G^\prime}{\partial y_i}$ under the composition
$$
 H^0(\P^4,\O_{\P^4}(2))\stackrel \rho\longrightarrow
 H^0(S,\O_{S}(2))\stackrel {g^*}\longrightarrow
 H^0(\P^2,\O_{\P^2}(6))
$$
where $\rho$ is restriction.  This composition is an injection
because $S$ is not contained in a quadric. So the $\frac{\partial
G^\prime}{\partial y_i}$ are also scalar multiples of each other,
say $\frac{\partial G^\prime}{\partial y_i}=\alpha_i
\frac{\partial G^\prime}{\partial y_0}$ for scalars $\alpha_i$. By
Euler's formula,
 $$
 3G^\prime \ \ =\ \
 \sum_{i=0}^4 y_i\  \frac{\partial G^\prime}{\partial y_i} \ \ =\ \
 \frac{\partial G^\prime}{\partial y_0}\  \sum_{i=0}^4 \alpha_iy_i
 $$
which implies $Y^\prime$ is not integral.  So $S=Y\cap Y^\prime$
is not integral.  This contradiction finishes the proof.
\end{proof}
\end{lem}

Now we make several observations, which we list in the following
Lemma.

\begin{lem}\ \\
(1) If $\Lambda\subset\P^4$ is any 2-plane, then $g^{-1}(\Lambda
\cap S)$ does not contain a curve of degree 3 or more.\\
(2) For $C^\prime\subset C$ any curve, the secant variety of
$C^\prime$ is contained in $Y^\prime$.\\
(3) Every curve in $C$ that lies in a hyperplane is a
plane curve.\\
(4) $C$ does not contain a conic curve.\\
(5) For every $x\in\P^2$, $g$ induces a nonzero map
$T_x\P^2\rightarrow T_{g(x)}Y$.
\end{lem}
\begin{proof}
(1) follows from the fact that the preimage of every hyperplane
containing $\Lambda$ is a distinct degree 3 curve in $\P^2$.  The
preimages are distinct because $S$ is not contained in a
hyperplane, which is our assumption for Case 1.
 \\
(2) holds because every line meeting $Y^\prime$ at two singular
points is contained in $Y^\prime$, and $Y^\prime$ is singular
along $C^\prime$.\\
(3) follows from (2) because $Y^\prime$ does not contain a
hyperplane.  This is because $\deg Y^\prime =3$ and $Y\cap
Y^\prime=S$ is integral.
 \\
(4) Suppose $C^\prime\subset C$ is a conic curve.  Then
$g^{-1}(C^\prime)$ contains a curve of degree at least 2 because
lines in $\P^2$ map to either lines or cubics. If $\Lambda$ is the
2-plane containing $C^\prime$, then $\Lambda\subset Y^\prime$ by
(2). So $\Lambda\cap S=\Lambda\cap Y$ is a degree 3 plane curve
consisting of $C^\prime$ and a line. This line is the image of a
curve in $\P^2$.  Therefore, $g^{-1}(\Lambda)$ contains a curve of
degree at least 3,
contradicting (1).\\
(5) If $T_x\P^2\rightarrow T_{g(x)}Y$ is the zero map, then every
line in $\P^2$ through $x$ maps to either a line or a cuspidal
plane cubic with cusp at $g(x)$.  In particular, the image of
every line in $\P^2$ through $x$ is contained in $T_{g(x)}Y$,
which is impossible by the assumption for this subsection that $S$
is not contained in a hyperplane.
\end{proof}

We will describe the integral curves $D$ that might occur in
$D(g)$ in terms of their multiplicity in $D(g)$ and in the
scheme-theoretic preimage $g^{-1} g(D)$.  For this purpose, we
recall the following definition.

\subsection*{Definition}
If $D$ is an integral curve contained in a one dimensional scheme
$Z$, then $mult_D Z$ will denote the length of the Artin ring
obtained by localizing $O_{Z}$ at the generic point of $D$.

\vspace{.35cm}

\begin{prop}
If $D\subset D(g)$ is an integral curve with $n:=\deg g|_D$, then
one of the following four cases holds:
\begin{itemize}
 \item[(1)] $n=1$,\ $\mult_D g^{-1} g(D)=2$,\ $\mult_DD(g)=1$,
    $D$ is a nonsingular conic, and $D$ is the only integral curve in $g^{-1} g(D)$.
 \item[(2)] $n=1$,\ $\mult_D g^{-1} g(D)=1$,\ $\mult_DD(g)=1$,
    and $g^{-1} g(D)$ contains two distinct nonsingular conics $D$, $D^\prime$.
 \item[(3)] $n=2$,\ $\mult_D g^{-1} g(D)=1$,\ $\mult_DD(g)=1$,
    $D$ is a nonsingular conic, and $g(D)$ is a singular plane cubic.
 \item[(4)] $n=3$,\ $\mult_D g^{-1} g(D)=1$,\ $\mult_DD(g)=2$,
    and $D$, $g(D)$ are both lines.
\end{itemize}
\end{prop}
\begin{proof}
We calculate $\mult_D D(g)$ as follows.  First suppose that $g$ is
unramified along every curve in $g^{-1} g(D)$ in the sense that
$T_x\P^2\rightarrow T_{g(x)}Y$ is injective for $x$ a general
point on any irreducible curve in $g^{-1} g(D)$. In this case, if
$D_1,\dots, D_r$ are the distinct integral curves contained in
$g^{-1} g(D)$, then
 \begin{equation}
 \mult_D D(g)=-1+\sum_{i=1}^r \deg g|_{D_i}
 \end{equation}
because the number of preimage points of a general point in $g(D)$
is the sum of the $\deg g|_{D_i}$. On the other hand, if $D$ is
the only integral curve contained in $g^{-1} g(D)$ and $g$ is
simply ramified along $D$ in the sense that $2D$ is contained in
$g^{-1} g(D)$ but $3D$ is not, then
 \begin{equation}
 \mult_D D(g)=1.
 \end{equation}
These are the only cases that will arise in our discussion here.
See \cite[Section 9.3]{F} for the construction of the Weil divisor
$D(g)$ representing $\D(g)$.

We claim that if $D\subset D(g)$ is a line, then $g(D)$ is a line
too. So suppose $D\subset D(g)$ is a line.  By Lemma 2.5(3),
$g(D)$ is not a twisted cubic. So it suffices to show that $g(D)$
is not a plane cubic. If it were, then by Lemma 2.5(2), $Y^\prime$
contains the 2-plane $\Lambda$ spanned by $g(D)$. Choose any
hyperplane $\Gamma$ containing $g(D)$, and let $\Sigma:=\Gamma\cap
Y^\prime$.  If $p\in\Gamma\cap Y^\prime$ is a singular point of
$Y^\prime$ away from $\Lambda$, then $\Gamma\cap Y^\prime$
contains the cone $\Sigma_p$ over $g(D)$ with vertex $p$.  Since
$\Lambda\cup\Sigma_p\subset \Sigma$ and $\Sigma$ and $\Sigma_p$
both have degree 3, no such $p$ can exist.  In particular,
$C\subset\Lambda$.  By Lemma 2.5(1), $g^{-1} (C)$ is either $D$,
$D$ with multiplicity 2, or $D$ and another line. So $D(g)$ will
have respectively degree 0, 1, or 2, according to the formulas for
$\mult_D D(g)$ above. However, $\deg D(g)=6$. So $g(D)$ is not a
plane cubic. Therefore, if $D\subset D(g)$ is a line, then $g(D)$
is too.

Now let $D\subset D(g)$ be any integral curve.  Using
$(n-1)D\subset D(g)$ and taking degrees, we obtain the inequality
$(n-1)\deg D\leq 6$.  The projection formula yields $3 \deg D = n
\deg g(D)$. Therefore, if $n\geq 2$, then
 \begin{equation}
 \frac{n}{3}\deg g(D)=\deg D\leq \frac{6}{n-1}.
 \end{equation}
This shows $n\leq 3$ because the left hand side is an integer. Now
we analyze what happens for each value $n=1, 2, 3$.

Suppose $n=3$, so that $\deg D=\deg g(D)\leq 3$.  If $\deg D=3$,
then Lemma 2.5(1) and (3) imply that $g(D)$ is neither a plane
cubic nor a twisted cubic.  So $\deg g(D)\leq 2$. By Lemma 2.5(4),
$\deg D\neq 2$.  So $D$, $g(D)$ are lines, as claimed in Case (4).

To see $\mult_D g^{-1} g(D)=1$, suppose $2D\subset g^{-1} g(D)$.
By Lemma 2.5(5), for $x\in D$ a general point, the map
$T_x\P^2\rightarrow T_{g(x)}Y$ has rank 1.  Choose $p\in g(D)$
general so that $(g|_D)^{-1}(p)$ consists of three distinct points
$x_1, x_2, x_3\in D$.  Through each $x_i$ there is a unique line
$l_i$ in $\P^2$ such that $T_{x_i}l_i\rightarrow T_pS$ is the zero
map. Each $g(l_i)$ is either a line through $p$ or a cuspidal
plane cubic with cusp at $p$.  In both cases, $g(l_i)$ is
contained in $T_pY$. Now the degree 3 curve $g^{-1}(T_pY)$
contains $2D+ l_1+ l_2+ l_3$. This contradiction proves $\mult_D
g^{-1} g(D)=1$.

To show $\mult_D D(g)=2$, it suffices by formula (2.1) to show
there is no integral curve $D^\prime\subset g^{-1} g(D)$ other
than $D$. If $D^\prime$ were such a curve, then $\deg
g|_{D^\prime}=3\deg D^\prime$ because $g(D)$ is a line. So for
$p\in g(D)$ a general point, $g^{-1}(p)$ would contain 3 points on
$D$ and $3\deg D^\prime$ points on $D^\prime$. There would be
$9\deg D^\prime$ lines joining one of the points on $D$ to one of
the points on $D^\prime$, and these lines would map to either
lines through $p$ or nodal plane cubics with node at $p$. These
image curves would be contained in $T_pY$, so we would have the
impossible situation that these $9\deg D^\prime$ lines would all
be contained in the degree 3 plane curve $g^{-1} (T_pY)$. So there
is no $D^\prime$, and $\mult_D D(g)=2$. This gives Case (4).

Suppose $n=2$.  If $2D\subset g^{-1} g(D)$, then $2D\subset g^{-1}
C$, which implies $2\deg D <6$ by Lemma 2.4.  By formula (2.3),
$g(D)$ can not be a line.  So $D$ is not a line, from above. Hence
$\deg D=2$ and $\deg g(D)=3$ by (2.3) and the fact $\deg D<3$.
Since $g(D)$ has degree 3 and can not be a twisted cubic, $g(D)$
is a singular plane cubic. By Lemma 2.5(1), $2D$ can not be
contained in $g^{-1} g(D)$, so $\mult_D g^{-1} g(D)=1$. Likewise,
$D$ is the only integral curve contained in $g^{-1} g(D)$, so
$\mult_D D(g)=1$ by formula (2.1). This gives Case (3).

Suppose $n=1$.  By (2.3), $\deg g(D)=3\deg D$. Since $g(D)$ is not
a line, $D$ is not a line.  So $\deg g(D)\geq 6$. If $D^\prime$ is
another integral curve contained in $g^{-1} g(D)$, then
$D^\prime\subset D(g)$.  By Cases (3) and (4), $\deg
g|_{D^\prime}=1$ because $\deg g(D^\prime)\geq 6$.  Therefore,
$\deg D^\prime=3\deg g(D)=\deg D$.

If $D$ is the only integral curve in $g^{-1} g(D)$, then
$2D\subset g^{-1}g(D)$ because $D\subset D(g)$ and $n=1$. By Lemma
2.4, $g^{-1}g(D)$ does not contain a curve of degree 6. So the
only possibility is $\mult_D g^{-1}g(D) = 2$ and $\deg D=2$
because $\deg D>1$. Since $D$ is integral,  $D$ is a nonsingular
conic. This gives Case (1) by formula (2.2).

The same reasoning shows that if $D^\prime$ is another integral
curve in $g^{-1} g(D)$, then the maximal curve contained in
$g^{-1}g(D)$ is $D+D^\prime$, and $D$, $D^\prime$ are both
nonsingular conics.  This gives Case (2) by formula (2.1).
\end{proof}

\begin{cor}
There are distinct lines $l_1, l_2, l_3$ in $\P^2$ such that
$$
D(g)=2l_1+2l_2+2l_3
$$
and the $g(l_i)$ are distinct lines in $S$.
\begin{proof}
In cases (1), (2), (3) of Proposition 2.6, $D$ has multiplicity 1
in $D(g)$. So if $D(g)$ does not contain a line, then $D(g)$ is
reduced and $D(g)\subset g^{-1}(C)$. This contradicts Lemma 2.4.
So $D$ contains a line $D$.  If $L=g(D)$ is the image line of $D$,
then $D$ is the only curve contained in $g^{-1}L$ by Case (4) of
Proposition 2.6.

Now suppose $D^\prime\neq D$ is another integral curve in $D(g)$.
We will show that $D^\prime$ is also a line. So suppose $\deg
D^\prime>1$.  Then $g(D^\prime)$ is not a line by Proposition 2.6.
Since $D$ and $D^\prime$ meet at a point, so do $L=g(D)$ and
$g(D^\prime)$.  If $g(D^\prime)$ is contained in a hyperplane,
then by Lemma 2.5(3), $g(D^\prime)$ is contained in a 2-plane. So
either $L$ and $g(D^\prime)$ are contained in the same 2-plane,
which is impossible by Lemma 2.5(1), or the curve $L\cup
g(D^\prime)$ spans a hyperplane, which is impossible by Lemma
2.5(3). So $g(D^\prime)$ is not contained in a hyperplane.

Let $p\in g(D^\prime)$ be a general point, and let $\Lambda_p$ be
the 2-plane spanned by $p$ and $L$.  Since the cubic $Y^\prime$ is singular
along $L$ and at $p$, $\Lambda_p$ is contained in
$Y^\prime$. Therefore, $\Lambda_p\cap S=\Lambda_p\cap Y$ is equal
to $L+Q_p$ for some conic curve $Q_p$ contained in $\Lambda_p$.

If $L$ is contained in $Q_p$, then $2L$ is contained in
$\Lambda_p\cap Y$.  So for every point $q\in L$,
$T_q\Lambda_p = T_q2L \subset T_qY$.  Since $g(D^\prime)$ is not contained
in any hyperplane, there is no hyperplane that contains
$\Lambda_p$ for every $p\in g(D^\prime)$.  Therefore, $\dim T_qY
=4$ for every $q\in L$ because $T_q\Lambda_p$ is
contained in $T_qY$ for $p\in g(D^\prime)$ a general point.  Since $Y$ is
nonsingular, one concludes that $L$ is not contained in $Q$.

If $Q_p$ is not a double line, then $g^{-1}Q_p$ contains a curve
of degree at least 2.  This is impossible by Lemma 2.5(1), because
$\Lambda_p$ contains $L$ and $Q_p$.  So for every $p\in
g(D^\prime)$, $Q_p=2L_p$ for some line $L_p\neq L$.

Therefore, $g(D^\prime)$ parametrizes a one dimensional family of
lines on $Y$.  By Proposition 2.6, every curve in $D(g)$ is
rational, so $g(D^\prime)$ is rational. However, the space of
lines on $Y$ does not contain a rational curve by Proposition 2.1.

This proves that every integral curve in $D(g)$ is a line, and
every line in $D(g)$ occurs with multiplicity 2 by Proposition
2.6.
\end{proof}
\end{cor}

\begin{prop}
The case $\deg Y^\prime=3$ does not occur.
\begin{proof}
We will use the $l_i$ from Corollary 2.7 to give an explicit
formula for $g$ and derive a contradiction using this formula.

Let $L_i:=g(l_i)$.  Each of the $L_i$ intersect the other two, but
they are not all contained in a plane by Lemma 2.5(1).  So they
all meet at some point $p$ in $S$, and $T_pY$ is the unique
hyperplane containing all the $L_i$.

I claim the $l_i$ do not all meet at a point $x$.  So suppose
there were such a point $x$, and note $g(x)=p$.  By Lemma 2.5(5),
at most one of the $g|_{l_i}$, say $g|_{l_1}$, ramifies at $x$.
Therefore there are points $x_2\in l_2$ and $x_3\in l_3$ different
from $x$ such that $g(x_2)=g(x_3)=p$, and $g$ maps the line
$l_{23}$ containing $x_2$ and $x_3$ to either a line or a nodal
plane cubic with node at $p$.  In particular, $g(l_{23})$ is
contained in $T_pY$.  Now the degree 3 curve $g^{-1}(T_pY)$
contains $l_1+ l_2+ l_3+ l_{23}$, which is impossible. So the
$l_i$ do not all meet at $x$.

Let $x_{ij}:=l_i\cap l_j$, and note $g(x_{ij})=p$.  If $g|_{l_1}$
is not ramified at either $x_{12}$ or $x_{13}$, then there is a
point $x^\prime$ on $l_1$ different from $x_{12}$, $x_{13}$ such
that $g(x^\prime)=p$, and the line from $x^\prime$ to $x_{23}$
gives a contradiction just as the line $l_{23}$ did above.  Hence
each $g|_{l_i}$ ramifies at exactly one of the $x_{ij}$.  So the scheme-theoretic preimage $g^{-1}(p)$ consists of three copies of $\Spec
k[\epsilon]/(\epsilon^2)$, one supported at each of the $x_{ij}$.  The scheme $g^{-1}(p)$ is
contained in the scheme $l_1+l_2+l_3$, but with no two copies of $\Spec
k[\epsilon]/(\epsilon^2)$ contained in the same $l_i$ because $(g|_{l_i})^{-1}(p)$ has length 3, as $g|_{l_i}$ is a morphism of degree 3 of nonsingular curves, and is therefore flat. Hence,
we can choose homogeneous coordinates $x_0, x_1, x_2$ on $\P^2$
such that $l_i=V(x_{i-1})$ and
$$
g^{-1}(p)=V(x_0^2x_1,\ x_1^2x_2,\ x_2^2x_0,\ x_0x_1x_2).
$$

There is a 3 dimensional space of hyperplanes in $\P^4$ containing
$p$, and these hyperplanes pull back to the linear system of
cubics on $\P^2$ spanned by
$$x_0^2x_1,\ x_1^2x_2,\ x_2^2x_0,\ x_0x_1x_2.$$
Set $p=(0,0,0,0,1)$ so that for a suitable choice of
coordinates $y_0,\dots, y_4$ on $\P^4$ we have
 \begin{equation}
 g=(x_0^2x_1,\ x_1^2x_2,\ x_2^2x_0,\ x_0x_1x_2,\ g_4)
 :\P^2\rightarrow \P^4
 \end{equation}
for some homogeneous polynomial $g_4$ of degree 3.

Since $g^{-1}V(y_3)=l_1+l_2+l_3$, one sees that $T_pY=V(y_3)$.  So
if $G$ is the homogeneous equation for $Y$, then
\begin{equation}
G=y_4^2y_3+y_4G_2+G_3
\end{equation}
where $G_i=G_i(y_0,\dots,y_3)$ is homogeneous of degree $i$.

By equation (2.4), the fact that $g$ is defined at $(1,0,0)\in
\P^2$ tells us that $x_0^3$ has nonzero coefficient in $g_4$.
Considering equation (2.4), we see that $x_0^7 x_1x_2$ has nonzero
coefficient in $g^*(y_4^2y_3)$.  Since $g(\P^2)\subset Y$, $g^*G$
is the zero polynomial, so $x_0^7x_1x_2$ has nonzero coefficient
in $g^*(y_4G_2+G_3)$ by equation (2.5).  The highest power of
$x_0$ that can appear in $g^*G_3$ is $x_0^6$ by (2.4). Therefore,
$x_0^7x_1x_2$ has nonzero coefficient in $g^*(y_4G_2)$.

The highest power of $x_0$ that can appear in a monomial of
$g^*G_2$ is $x_0^4$, and this necessarily occurs in the monomial
$g_0^2=x_0^4x_1^2$.  So if $x_0^7$ appears in a monomial of
$g^*(y_4G_2)=g_4 g^*G_2$, then that monomial is $x_0^7x_1^2$, not
$x_0^7x_1x_2$.  This contradiction shows that $\deg
Y^\prime\neq3$.
\end{proof}
\end{prop}

\subsection*{Case 2}
By Proposition 2.8,  the image of $g:\P^2\rightarrow Y$ is the
surface $S:=Y\cap Y^\prime$, where $Y^\prime$ is a hyperplane in
$\P^4$.  We will derive a contradiction, thus proving Theorem 2.

The morphism $\delta: Y\rightarrow Y^*$ sending a point $p$ to the
hyperplane $T_pY$ is given by an ample invertible sheaf on $Y$,
since $\O_Y(1)$ generates $\Pic Y$ by the Grothendieck-Lefschetz
Theorem, \cite[Theorem 4.3.2]{H1}. So $\delta$ is a finite
morphism. In other words, no hyperplane in $\P^4$ can be tangent
to $Y$ at infinitely many points. Therefore $S$ is a cubic surface
with isolated singularities.

For $s\in S$ a general closed point, let $\{E_t\}$ be a general
pencil of hyperplane sections of $S$ containing $s$, and let
$E'_t:=g^{-1} (E_t)$. Then $E_t$ and $E^\prime_t$ are both smooth
plane cubics for $t\in \P^1$ general.

Since $g:\P^2\rightarrow S$ has degree 3, we can set
$g^{-1}(s)=\{a, b, c\}$.  Then $(E^\prime_t, a)\rightarrow (E_t,
s)$ is an isogeny of elliptic curves of degree 3 for $t$ general.
Note that $b$ is a 3-torsion point of $E_t'$.

Let $E^\prime\subset\P^1\times\P^2$ and $E\subset\P^1\times S$ be
the total spaces of the families $\{E_t^\prime\}$ and $\{E_t\}$
over the base $\P^1$, and consider the diagram
$$
\begin{array}{ccc}
 E^\prime & \longrightarrow & \P^2      \\
 \downarrow\rlap{$g^\prime$}  & & \downarrow\rlap{$g$}             \\
 E & \longrightarrow & S
\end{array}
$$
where $g^\prime$ is the restriction of $\text{id}_{\P^1}\times g$
to $E^\prime\subset \P^1\times\P^2$.

Let $\zeta$ be the generic point of $\P^1$, and let
$E^\prime_\zeta$, $E_\zeta$ be the fibers of $E^\prime$, $E$ over
$\zeta\in\P^1$. Then $A:=\zeta\times a$ and $B:=\zeta\times b$ are
$\zeta$-valued points of $E^\prime_\zeta$, and $B$ is a 3-torsion
point of the elliptic curve $(E^\prime_\zeta, A)$. Hence $P\mapsto
P+B$ gives an automorphism of $E^\prime_\zeta$ over $E_\zeta$ of
order 3 that extends to some birational self-map $\psi:
E^\prime\dashrightarrow E^\prime$ of $E^\prime$ over $E$. In other
words, $g^\prime=g^\prime\circ\psi$ as rational maps. The
horizontal arrows in the diagram above are birational morphisms.
So $\psi$ induces an order 3 birational self-map
$\phi:\P^2\dashrightarrow\P^2$ of $\P^2$ such that $g=g\circ\phi$
as rational maps.  Therefore $\phi$ is an automorphism of $\P^2$
by the following lemma.

\begin{lem}
If $\phi:\P^2\dashrightarrow\P^2$ is a birational map and
$g:\P^2\rightarrow S$ is a dominant morphism to a surface $S$ such
that $g\circ \phi=g$ as rational maps, then $\phi$ extends to an
automorphism of $\P^2$.
\end{lem}
\begin{proof}
Let $U\subset\P^2$ be the domain of definition of $\phi$.  Then
$U$ is the compliment of finitely many points in $\P^2$, and $\Pic
U$ is generated by $\O_U(1)$, the restriction of $\O_{\P^2}(1)$ to
$U$. From $g\circ \phi=g$ it follows that
$\phi^*g^*\O_S(1)=g^*\O_S(1)|_U$.  In other words, $\phi^*
\O_{\P^2}(3)= \O_U(3)$.  Therefore $\phi^*\O_{\P^2}(1)=\O_{U}(1)$.

Since $\phi$ is dominant, $\phi:U\rightarrow \P^2$ is given by
three linearly independent sections $\sigma_0, \sigma_1, \sigma_2$
of $\O_U(1)$, and these sections are the restriction to $U$ of
sections $\tau_0, \tau_1\, \tau_2$ of $\O_{\P^2}(1)$.  The
$\tau_i$ are linearly independent, and hence there is no point
where they all vanish.  So $\phi$ extends to the automorphism
$(\tau_0, \tau_1, \tau_2)$ of $\P^2$.
\end{proof}

Since $\phi$ has order 3, its matrix is one of the following after
being placed in Jordan form and scaling:
\begin{gather*}
 \phi=
\begin{bmatrix}
    1   &       &       \\
        & \rho  &       \\
        &       & \rho
\end{bmatrix}
 \quad \text{or} \quad
 \phi=
\begin{bmatrix}
    1   &       &       \\
        & \rho  &       \\
        &       & \rho^2
\end{bmatrix}
\end{gather*}
where $\rho^3=1$, $\rho\neq 1$.  Fix an identification of
$Y^\prime$ with $\P^3$ so that we can write $g:\P^2\rightarrow
S\subset Y^\prime$ as $(g_0,\dots,g_3)$ for some degree 3
polynomials $g_i$. Then
$(g_0,\dots,g_3)=(\phi^*g_0,\dots,\phi^*g_3)$. So each of the
$g_i$ are eigenvectors of $\phi^*:H^0(\P^2,\O(3))\rightarrow
H^0(\P^2,\O(3))$ with the same eigenvalue.

Suppose
$$
\phi =
 \left[
  \begin{array}{ccc}
    1   &       &       \\
        & \rho  &       \\
        &       & \rho
  \end{array}
 \right]
$$
\smallskip\smallskip
Then the following are a basis for each of the eigenspaces of
$\phi^*$:
\begin{align*}
1 :         & \ \ x_0^3,\ x_1^3,\ x_2^3,\ x_1^2x_2,\ x_1x_2^2    \\
\rho :      & \ \ x_0^2x_1,\ x_0^2x_2      \\
\rho^2 :    & \ \ x_0x_1^2,\ x_0x_2^2,\ x_0x_1x_2
\end{align*}
Since $S$ is not a 2-plane in $\P^3$, the eigenspace containing
the $g_i$ has dimension at least 4. So the $g_i$ all have
eigenvalue 1.

Now consider the morphism
 $$
 h=(x_0^3,\ x_1^3,\ x_2^3,\ x_1^2x_2,\ x_1x_2^2):\P^2\rightarrow
 \P^4.
 $$
The image surface $S^\prime$ of $h$ is a cone over a twisted cubic
curve $T^\prime$ with vertex $(1,0,0,0,0)$. The morphism
$g:\P^2\rightarrow Y^\prime\cong\P^3$ is the morphism $h$ followed
by projection $\pi_p$ from some point $p$. Projection maps lines
to lines, so the image $S$ of $g$ is a cone over some cubic plane
curve $T$. The curve $T$ is singular because it is the image of a
rational twisted cubic curve. Therefore, $S$ is singular along a
line. We already saw that $S$ has only finitely many
singularities. This contradiction rules out the first possibility
for $\phi$.

Therefore
$$
\phi =
 \left[
  \begin{array}{ccc}
    1   &       &       \\
        & \rho  &       \\
        &       & \rho^2
  \end{array}
 \right]
$$
\smallskip\smallskip
So the following are bases for the eigenspaces of $\phi^*$:
\begin{align*}
1 :         & \ \ x_0^3,\ x_1^3,\ x_2^3,\ x_0x_1x_2    \\
\rho :      & \ \ x_0^2x_1,\ x_0x_2^2,\ x_1^2x_2      \\
\rho^2 :    & \ \ x_0^2x_2,\ x_0x_1^2,\ x_1x_2^2
\end{align*}
The eigenspace containing the $g_i$ has dimension at least 4.  So
$g:\P^2\rightarrow Y^\prime=\P^3$ is the morphism $(x_0^3,\
x_1^3,\ x_2^3,\ x_0x_1x_2)$.  Embed $Y^\prime$ in $\P^4$ as
$V(y_4)$ so that $g:\P^2\rightarrow \P^4$ has the form claimed in
Theorem 2.

\section{Morphisms from Quintic to Cubic Threefolds}

In Section 3.1 we discuss preimages of lines on $Y$ under the
rational map $H\dashrightarrow Y$. Section 3.2 gives some
information about $H$ in case it has a component that does not map
dominantly onto $Y$. In Section 3.3, we consider the various
possibilities for $H$ and rule them out case by case. When $H$ is
integral, the results of Section 3.1 will be the main tool. When
$H$ is more degenerate, Theorem 2 will play a central role.

\subsection{Preimages of Lines on $Y$}

First we will prove a basic fact about how the dualizing sheaf of
a curve behaves under normalization.  Then we will discuss the
family of lines on $Y$ and what can be said about the preimage in
$H$ of a general line in $Y$.  The main result is Corollary 3.6.

\begin{lem}
Let $C$ be a projective Gorenstein scheme of pure dimension 1 that
is reduced at the generic point of one of its irreducible
components $C^\prime$.  If $\nu:\tilde{C}\rightarrow C^\prime$ is
the normalization map of the reduced structure on $C^\prime$ and
$N$ is the number of points $p\in\tilde C$ such that the map of
local rings $\O_{C,\nu(p)}\rightarrow\O_{\tilde C,p}$ fails to be
an isomorphism, then
$$
\deg \nu^*\omega^o_C\ \geq\  \deg \omega_{\tilde{C}} + N.
$$
If $N=0$, then equality holds.
\begin{proof}
Following \cite[Ex. 3.7.2]{H2}, we compute
\begin{align*}
\omega_{\tilde C} &\cong \Hom_C(\nu_*\O_{\tilde C}, \omega_C^o) \\
& \cong \Hom_C(\nu_*\O_{\tilde C}, \O_C)\otimes_{\O_C} \omega_C^o
\end{align*}
as $\O_{\tilde C}$-modules.  Since $\omega_{\tilde C}$ and
$\omega_C^o$ are invertible, $\Hom_C(\nu_*\O_{\tilde C}, \O_C)$ is
an invertible sheaf of $\O_{\tilde C}$-modules.  So it is enough
to show that $\Hom_C(\nu_*\O_{\tilde C}, \O_C)$ is isomorphic to
an ideal sheaf in $\O_{\tilde C}$ corresponding to a closed
subscheme of $\tilde C$ supported at the points where
$\nu:\tilde{C}\rightarrow C$ fails to be an isomorphism.

Consider the map of local rings
 $$A:=\O_{C,\nu(p)} \longrightarrow B:=\O_{\tilde C,p}$$
for some closed point $p\in \tilde C$. We will show that
$\Hom_A(B,A)$ is an ideal in $B$ and is the unit ideal if and only
if $A\rightarrow B$ is an isomorphism.

Note that $B$ is the normalization of $A/P$ in its field of
fractions, where $P\subset A$ is the prime ideal of
$C^\prime\subset C$. Consider the map of $B$-modules
 $$\Phi:\Hom_A(B,A)\longrightarrow B$$
that sends $\psi\in \Hom_A(B,A)$ to the equivalence class
$\overline{\psi(1)}$ of $\psi(1)$ in $A/P$, which injects into
$B$.

If $\Phi(\psi)=0$, then $\psi(1)\in P$, whence $\psi(B)$ is
contained in the ideal $P\subset A$.  In other words, $\psi$ is a
local section of the $\O_C$-module $\Hom_C(\nu^*\O_{\tilde
C},\O_C)$ defined in a neighborhood of $\nu(p)$ in $ C$ such that
$\psi$ vanishes on the reduced structure of $C^\prime$.  If we
consider $\Hom_C(\nu^*\O_{\tilde C},\O_C)$ as an $\O_{\tilde
C}$-module, then $\psi$ is a local section defined on a
neighborhood of $p$ in $\tilde C$, and $\psi$ vanishes at all but
finitely many points in the neighborhood. Since
$\Hom_C(\nu^*\O_{\tilde C},\O_C)$ is an invertible $\O_{\tilde
C}$-module, $\psi=0$.  Therefore $\Phi$ is an injection of
$B$-modules, whence $\Hom_A(B,A)$ is realized as an ideal in $B$.
It remains to check that $A\rightarrow B$ is an isomorphism if and
only if $\Phi$ is.

If $A\rightarrow B$ is an isomorphism, then $\Phi$ is clearly an
isomorphism. Conversely, suppose $\Phi$ is an isomorphism.  Then
there is some $\psi\in \Hom_A(B,A)$ such that $\Phi(\psi)=1$, i.e.
$\psi(1)-1\in P$. Therefore $\psi(1)$ is a unit in $A$ because $A$
is local. Hence $\Psi : \Hom_A(B,A)\rightarrow A$ given by
$\psi\mapsto\psi(1)$ is surjective because it is a morphism of
$A$-modules.  Note that $\Phi$ factors as $\Psi$ followed by
$A\rightarrow B$.  Since $\Psi$ is surjective and $\Phi$ is
injective, we conclude that $A\rightarrow B$ is injective.  And
since $\Phi$ is surjective, $A\rightarrow B$ is surjective.
\end{proof}
\end{lem}

\begin{lem}
There are only finitely many closed points $p\in Y$ such that
there are infinitely many lines on $Y$ through $p$.
\begin{proof}
Suppose not.  Then there is a curve $C\subset Y$ such that for
every point $p\in C$ there are infinitely many lines on $Y$
through $p$.  So for every $p\in C$ there is an irreducible
component $\Sigma_p$ of $Y\cap T_pY$ such that $\Sigma_p$ is a
cone over a plane curve with vertex $p$.

Recall from \cite {HRS} or \cite{S2} that the family $\Delta_1$ of
lines on $Y$ is an irreducible surface, so a dimension count shows
that every line on $Y$ lies on one of the surfaces $\Sigma_p$ for
some $p\in C$. Also recall from \cite{S2} that a general line $L$
in $Y$ has normal bundle $N_{L/Y}=\O_L\oplus\O_L$.

Fix a general line $L$ on $Y$, lying on $\Sigma_p$.  There is some
2-plane $\Lambda$ that is tangent to $\Sigma_p$ at every point of
$L$ on account of $\Sigma_p$ being a cone. Therefore
$T_q\Lambda\subset T_qY$ for every $q\in L$. It follows from
Nakayama's Lemma that these pointwise inclusions give an injection
$T_{\Lambda}|_L\rightarrow T_Y|_L$ of $\O_L$-modules. Note that
$T_{\Lambda}|_L=\O_L(1)\oplus\O_L(2)$, and consider the normal
bundle sequence
$$
0\longrightarrow T_L\longrightarrow T_Y|_L
\stackrel{\phi}{\longrightarrow} \O_L\oplus\O_L\longrightarrow 0.
$$
By the above description of $T_\Lambda|_L$, the composition
 $$T_\Lambda|_L\longrightarrow T_Y|_L\stackrel\phi\longrightarrow \O_L\oplus\O_L$$
is the zero morphism.  So $\ker\phi$ has rank at least 2,
contradicting $\ker\phi=T_L$.
\end{proof}
\end{lem}

We will need the following modification of \cite[Lemma 2.1]{S2}.

\begin{lem}
Let $g:Z\rightarrow Y$ be a morphism from a purely 3 dimensional
separated scheme of finite type over the ground field. If $L$ is a
general line on $Y$, then $g^{-1}(L)$ has pure dimension 1.

Moreover, if $Z$ is integral and $g$ is dominant, then $g^{-1}(L)$
is singular at only finitely many points.  In other words,
$g^{-1}(L)$ is reduced at the generic point of each irreducible
component.
\begin{proof}
Let $\mathbb F$ be the total space of the family of lines on $Y$
with base $\Delta$. From \cite{S2}, $\Delta$ is a smooth surface,
and there is some open subscheme $\Delta_0$ in $\Delta$ with
preimage $\mathbb F_0$ in $\mathbb F$ such that $\mathbb F_0$ is a
locally trivial fiber bundle over $\Delta_0$, whose fibers are
lines in $Y$.

Let $Z_0$ be the union of the surfaces in $Z$ that are mapped to
points in $Y$, let $Z_1$ be the union of the curves on $Z$ that
are mapped to points in $Y$, and let $Z_2$ be the set of points in
$Z$ at which $g$ fails to induce an injection on tangent spaces.
Note that $Z_2$ contains the singular locus of $Z$.

Note that $g(Z_0)$ has dimension at most 1, and $g(Z_1)$ has
dimension at most 2. Since the canonical morphism $\mathbb
F_0\rightarrow Y$ is dominant, a general line will meet $g(Z_1)$
in only finitely many points. Also, if a general line meets
$g(Z_0)$, then $g(Z_0)$ is a curve such that for every point $p\in
g(Z_0)$ there are infinitely many lines on $Y$ through $p$. This
contradicts Lemma 3.2. So a general line does not meet $g(Z_0)$.
Therefore, $g^{-1}(L)$ has pure dimension 1 for $L$ a general line
on $Y$.

Now suppose that $Z$ is integral and $g$ is dominant. Then $\dim
g(Z_i)\leq i$ for $i=0, 1, 2$.  We use the characteristic zero
assumption of this section to get $\dim g(Z_2)\leq 2$. So a
general line will meet $g(Z_2)$ in only finitely many points, and
will not meet $g(Z_0)$ and $g(Z_1)$ at all. Hence $g^{-1}(L)$ is
nonsingular away from the preimage of $L\cap g(Z_2)$, which
consists of only finitely many points.
\end{proof}
\end{lem}

\subsection*{Definition}
If $Z$ is any scheme and $F:Z\dashrightarrow \P^n$ is a rational
map given by sections $F_0, \dots, F_n$ of some line bundle on
$Z$, then let \emph{indet$(F)$} denote the scheme of common
vanishing of the $F_i$ in $Z$:
$$indet(F):= V(F_0,\dots, F_n)\subset Z.$$

\begin{lem}
Take $F$ is as in the definition, and let $\pi:\widetilde
Z\rightarrow Z$ be the blowup of $Z$ in $\indet (F)$.  Then there
is a canonical morphism $F:\widetilde Z\rightarrow \P^n$ such that
$\widetilde F=F\circ\pi$ as rational maps.  Moreover, if
$p\in\indet(F)$ is a closed point, then $\widetilde F$ induces a
closed immersion of $\pi^{-1}(p)$ into $\P^n$.
\end{lem}
\begin{proof}
Recall that $\widetilde Z$ is isomorphic to the closure of the
graph of $F$, $\Gamma_F\subset Z\times \P^n$. Thus projection onto
$\P^n$ induces the desired morphism $\widetilde F:\widetilde
Z\rightarrow \P^n$. If $p\in \indet(F)$ is a closed point, then
$\pi^{-1}(p)$ is a closed subscheme of $p\times\P^n$, which maps
isomorphically onto $\P^n$ by projection.
\end{proof}

\begin{lem}
If $p\in\indet(F)$ is a closed point, then $p\in H$ is a point of
order at least 3.
\end{lem}
\begin{proof}
Let $F=(F_0,\dots,F_4)$, and let $Y$ have homogeneous equation
$G$.  Since the $F_i$ vanish at $p$ and $G$ has degree 3,
$G(F_0,\dots, F_4)$ vanishes to order at least 3 at $p$.  In other
words, $p\in F^{-1}(Y)$ is a point of order at least 3.  But $p$
is contained in $H$ and not $X$.  The lemma follows.
\end{proof}

If $L$ is a general line in $Y$ with $L:=V(\xi_1, \xi_2, \xi_3)$
for some linear forms $\xi$ on $\P^4$, then define
\begin{align*}
 F^{-1}(L) & :=V(F^*\xi_1, F^*\xi_2,F^*\xi_3) \\
 C & :=F^{-1}(L)\cap H \\
 D & :=F^{-1}(L)\cap X
\end{align*}
If $\widetilde H\rightarrow H$ is the blowup of $H$ at the
indeterminacy scheme $\indet(F|_H)$ and $h:\widetilde H\rightarrow
Y$ is the resulting morphism, then Lemma 3.3 says that $h^{-1}(L)$
is purely one dimensional, so the same holds for $C$. Lemma 3.3
implies $D$ is also purely one dimensional. So $F^{-1}(L)$ is a
complete intersection in $\P^4$.

\begin{lem}
With $C$ as above,  $\omega_C^o=\O_C(-1)$.
\begin{proof}
Since $F^{-1}(L)$ is the complete intersection of three cubic
hypersurfaces in $\P^4$,
\begin{align*}
 \omega^o_{F^{-1}(L)} &=\O_{F^{-1}(L)}(-5+3+3+3) \\
 &=\O_{F^{-1}(L)}(4).
\end{align*}
Following \cite[Ex. 3.7.2]{H2}, compute
\begin{align*}
\omega^o_C &= \Hom_{F^{-1}(L)}(\O_C, \omega_{F^{-1}(L)}^o) \\
& = \Hom_{F^{-1}(L)}(\O_C, \O_{F^{-1}(L)})\otimes
\O_{F^{-1}(L)}(4).
\end{align*}
Therefore, if $I_C, I_D\subset\O_{F^{-1}(L)}$ denote the ideals of
$C, D\subset {F^{-1}(L)}$, then it suffices to carry out the
following computation
\begin{align}
 \Hom_{F^{-1}(L)}(\O_C, \O_{F^{-1}(L)})
 & \cong (0:I_C) \\
 & \cong I_D \\
 & \cong \O_C(-C\cap D) \\
 & \cong \O_C(-5).
\end{align}

The isomorphism (3.1) is given by $\psi\mapsto \psi(1)$.

To see (3.2), note that $D$ is Cohen-Macaulay since it is a local
complete intersection in $X$ and $X$ is C.M., as $L$ is a l.c.i.
in $Y$. It follows that $C$ and $D$ are linked because
$F\inverse(L)$ is Gorenstein, cf. \cite[Theorem 21.23]{E}.

To see (3.3), we will show that
$$
I_D\longrightarrow\frac{I_D+I_C}{I_C}
$$
is an isomorphism of $\O_C$-modules, i.e. $I_D\cap I_C=0$.  So let
$a$ be a local section of $I_D\cap I_C$, which is necessarily
supported on $C\cap D$.  We will show $a=0$.

Let $A=\O_{F^{-1}(L),p}$ for some closed point $p\in C\cap D$ with
maximal ideal $m$. Let $b$ denote the image of $a$ in $A$. Some
power of $m$ annihilates $b$ because $b$ is supported at $p$.
Therefore $m$ is the annihilator of some nonzero multiple of $b$,
provided $b\neq 0$.  So $m$ is an associated prime ideal of $A$.
But $A$ is C.M. and one dimensional, so every associated prime is
minimal.  Hence $b=0$ after all.  Therefore $a=0$ because its
localization is zero at every point.

The isomorphism (3.4) follows from $C\cap D=C\cap X$.
\end{proof}
\end{lem}

\begin{cor}
If $C:=F^{-1}(L)\cap H$ is reduced at the generic point of an
irreducible component $C^\prime$, then $C^\prime$ is either a
smooth plane conic disjoint from the rest of $C$, or $C^\prime$ is
a line meeting the rest of $C$ at only one point.
\begin{proof}
This is immediate from Lemma 3.1 and Lemma 3.5 because the degree
of the dualizing sheaf of a smooth curve is at least $-2$.
\end{proof}
\end{cor}

\subsection{A Multiplicity  Result}
The following result holds in arbitrary characteristic and will be
used in the next subsection.

\begin{lem}
Let $Z$ be the reduced structure on an irreducible component of
$H$.  Assume $m\leq d$.  If $F|_Z:Z\dashrightarrow Y$ is not
dominant, then $2Z\subset H$ as divisors in $\Pn$.
\begin{proof}
Suppose $F$ does not map $Z$ dominantly onto $Y$.  Then $Z$ is
covered by curves that are mapped to points under $F$.  So suppose
that $C\subset Z$ is an integral curve with $F(C)=(1,0,\dots,0)$
for simplicity, and let $I\subset\O_{\Pn}$ be the ideal sheaf of
$C$.  If $F=(F_0,\dots, F_n)$, then $F_i\in H^0(\Pn, I(m))$ for
$i>0$. We show that $H$ is singular along $C$.

Let $K\in H^0(\Pn, I(em-d))$ be the homogeneous equation for $H$.
Let $I^{(2)}\subset\O_{\Pn}$ be the ideal such that
$I/I^{(2)}=(I/I^2)/\text{torsion}$.  Then $K$ induces a section
$\overline{K}\in H^0(\Pn, I/I^{(2)}(em-d))$.  We will show
$\overline{K}=0$.

If $G=G(y_0,\dots,y_n)$ is the homogeneous equation of $Y$, then
we can write $G=y_0^{e-1}G_1+\dots+y_0G_{e-1}+G_e$ where the
$G_i=G_i(y_1,\dots,y_n)$ are homogeneous of degree $i$. So the
homogeneous equation of $F^{-1}(Y)=X+H$ is
\begin{equation}
F^*G=F_0^{e-1}F^*G_1+F_0^{e-2}F^*G_2+\dots
\end{equation}
where $F^*G=G(F_0,\dots,F_n)$.

Let $D$ be the Cartier divisor $V(F_0^{e-1})\cap C$ on $C$.  Note
that $D$ is disjoint from $X$ because $D$ is supported on $C\cap
V(F_0,\dots,F_n)$ and $V(F_0,\dots,F_n)$ is disjoint from $X$.  So
from (3.5) it follows that $\overline K$ restricts to the zero
section on $D$ because $F_0^{e-i}F^*G_i$ restricts to zero on $C$
for $i\geq 2$, and $F_0^{e-1}$ restricts to zero on $D$.

Tensor the exact sequence
$$
0\longrightarrow\O_C(-D)\longrightarrow\O_C\longrightarrow
\O_D\longrightarrow 0
$$
with $I/I^{(2)}(em-d)$ and use $\O_C(-D)=\O_C(-m(e-1))$ to obtain
$$
I/I^{(2)}(m-d)\stackrel{\tau}{\longrightarrow} I/I^{(2)}(em-d)
\stackrel{\rho}{\longrightarrow} I/I^{(2)}\otimes\O_D (em-d)
\longrightarrow 0
$$
To see that $\tau$ is an injection, note that $\tau$ is
multiplication by $F_0^{e-1}$, which is a unit in the local ring
of almost all the points of $C$.  Therefore $\tau$ could only have
torsion elements in its kernel.  But $I/I^{(2)}(em-d)$ is a
torsion-free sheaf.  So $\tau$ is injective.

Since $\rho(\overline{K})=0$, $\overline{K}$ is the image of a
section
 $$
 \widetilde{K}\in H^0(C,I/I^{(2)}(m-d)).
 $$
From the conormal sequence of $C$ in $\Pn$, we get a morphism
$I/I^{(2)}(m-d) \rightarrow \Omega^1_{\Pn}|_C(m-d)$ that is an
injection on the regular locus of $C$, and is therefore an
injection because $I/I^{(2)}(m-d)$ is torsion-free.

The sheaf $\Omega^1_{\Pn}|_C(m-d)$ has no nonzero global sections
because of the injection $\Omega^1_{\Pn}|_C(m-d)\rightarrow
\O_{\Pn}(m-d-1)^{\oplus n+1}$ from the Euler sequence. Therefore,
$I/I^{(2)}(m-d)$ has no nonzero global sections. So
$\widetilde{K}=0$, whence $\overline{K}=0$. Therefore $K\in
H^0(\Pn, I^{(2)})$, which implies $H$ is singular along $C$. Since
$Z$ is covered by such curves, $H$ is singular at every point of
$Z$. This is only possible if $2Z\subset H$ as divisors.
\end{proof}
\end{lem}

\subsection{Morphisms from Quintic to Cubic Threefolds}

Recall the decomposition $F^{-1}(Y)=X+H$ from Section 1. We will
consider the various possibilities for $H$ and rule them out one
at a time, thus proving Theorem 1.

\begin{prop}
$H$ does not contain a hyperplane in $\P^4$ that maps dominantly
onto $Y$.
\begin{proof}
Suppose $Z\subset H$ is a hyperplane that maps dominantly onto
$Y$. If we restrict $F$ to any $\P^2$ contained in $Z$, then
$F|_{\P^2}$ is described by Theorem 2.  In particular, $F(\P^2)$
is the intersection of $Y$ with a hyperplane tangent to $Y$ at 3
points.  The family of 2-planes contained in $Z$ and the family of
tangent planes to $Y$ both have dimension 3, so a general tangent
plane to $Y$ is tangent to $Y$ at 3 points. However, since $Y$ is
a nonsingular hypersurface, a general tangent plane to $Y$ is
tangent at only one point, cf. \cite[Lemma 5.15]{CG}.
\end{proof}
\end{prop}

\begin{cor}
If $K$ is a hyperplane contained in $H$, then $2K\subset H$, and
we can choose coordinates $x_0, \dots, x_4$ on the $\P^4$
containing $K$ such that $K=V(x_4)$ and $F|_K$ is given by the
formula
$$ F|_K= (x_0^3,\
x_1^3,\ x_2^3,\ x_0x_1x_2,\ 0)
$$
\begin{proof}
By Lemma 3.7 and Proposition 3.8, $2K\subset H$.

Since $F|_K$ is not dominant and $K\cong\P^3$, $F|_K$ can not be a
morphism.  Choose a point $p$ in the indeterminacy locus of
$F|_K$. Let $\pi:\widetilde{K}\rightarrow K$ be the blowup of $K$
at the indeterminacy scheme $\indet(F|_K)$, so that the rational
map $F|_K$ extends to a morphism $\Phi:\widetilde{K}\rightarrow
S$, where $S$ is the image surface of $K$ under $F$.  By Lemma
3.4, $\pi^{-1}(p)$ maps isomorphically onto $S$ under the morphism
$\Phi$.  So for $s\in S$ a general point, $\Phi^{-1}(s)$ is a
curve in $\tilde{K}$ that meets $\pi^{-1}(p)$. Therefore the
preimage of $s$ under $F|_K$ is a curve in $K$ through $p$.

If $\Lambda$ is any 2-plane in $K$ disjoint from $\indet(F|_K)$,
then $F|_\Lambda$ is described by Theorem 2.  In particular, for
$s\in S$ any nonsingular point, the preimage of $s$ in $\Lambda$
consists of three reduced points.  Therefore the preimage of $s$
in $K$ consists of three distinct lines that are reduced away from
$\indet (F|_K)$.  Indeed, if the preimage of $s$ in $K$ contained
a curve other than a line, we could choose $\Lambda$ such that
$F|_\Lambda^{-1}(s)$ was not 3 reduced points by choosing
$\Lambda$ to be tangent to $F|_\Lambda^{-1}(s)$ at a point but not
contain any component of $F|_\Lambda^{-1}(s)$.

So there is a two parameter family of lines in $K$ that are mapped
by $F$ to points in $S$. Every such line meets $\indet (F|_K)$,
which consists of finitely many points.  So there is some point
$p\in\indet(F|_K)$ such that there is a two parameter family of
lines in $K$ through $p$ that are each mapped to a point in $S$
under $F$.

A line $L$ is mapped to a point by $F$ exactly when the scheme
$L\cap\indet(F)$ has length 3 because $F^*\O(1)=\O(3)$.  Since a
general line in $\P^4$ through $p$ meets $\indet(F)$ in a scheme
of length 3, the same holds for every line through $p$.  So every
line through $p$ is mapped to a point.  Therefore, $F|_K$ is
determined by $F|_\Lambda$ for any 2-plane $\Lambda$ in $K$ not
containing $p$.  Theorem 2 determines $F|_\Lambda$.  If one takes
$p=V(x_0, x_1, x_2)$ and $\Lambda=V(x_3)$, where $x_0,\dots, x_3$
are homogeneous coordinates on $K$, then $F|_K$ has the desired
form.
\end{proof}
\end{cor}

\begin{prop}
$H$ is not equal to $Q+2K$ for an integral quadric $Q$ and a
hyperplane $K$.
\begin{proof}
Suppose $H=Q+2K$.  By Corollary 3.9, $F|_K$ factors through
projection from some point $p\in K$.  By Lemma 3.5, $H$ has order
at least 3 at every point in $\indet(F)$. So every point in
$\indet(F)$ is in $K\cap Q$. Hence, $p$ is the only point in
$\indet(F)$ because it is the only point on $K$ where $F$ is
undefined.

The map $F|_Q$ is dominant by Lemma 3.7.  Let $\pi:\widetilde
Q\rightarrow Q$ be the blowup of $Q$ at the indeterminacy scheme
$\indet(F|_Q)$, and let $q:\widetilde Q\rightarrow Y$ be the
resulting morphism extending $F|_Q$.  Apply Lemma 3.3 to $q$ to
see that for $L$ a general line, $q^{-1}(L)$ is reduced at the
generic point of each of its irreducible components, whence the
same is true for $C_Q:= F^{-1}(L)\cap Q$.

By Lemma 3.4, $q$ is a closed immersion when restricted to
$\pi^{-1}(p)$. Since $q(\pi^{-1}(p))$ is an effective divisor on
$Y$, it is ample by the Grothendieck-Lefschetz
Theorem,\cite[Theorem 4.3.2]{H1}. So every line in $Y$ meets
$q(\pi^{-1}(p))$.  Therefore, $C_Q$ has an irreducible component
$C_1$ containing the point $p$. By the same argument,
$F^{-1}(L)\cap K$ has an irreducible component $C_2$ that contains
$p$.  Since $L$ is a general line on $Y$, which is covered by
lines, $L$ is not contained in the image of $K\cap Q$. So $C_1\neq
C_2$.  By Lemma 3.6, $\omega_C^o=\O_C(-1)$.  Since $C$ has more
than one irreducible component than contains $p$, Lemma 3.1 shows
that every irreducible component of $C_Q$ that contains $p$ is a
line because the dualizing sheaf of every smooth curve is at least
$-2$.  So $Q$ is covered by lines through $p$, and these lines are
parametrized  by a general hyperplane section $Q\cap \P^3$ of $Q$.
Because the line $L\subset Y$ is general, $\dim Q\cap\P^3=2$, and
$\dim \Delta =2$, where $\Delta$ is the space of lines on $Y$, we
conclude that every general line on $Q$ through $p$ maps to a line
on $Y$.  Hence the rational map $F|_Q:Q\dashrightarrow Y$ induces
a rational map $Q\cap \P^3\dashrightarrow \Delta$.  However,
$\Delta$ does not contain a rational curve by Proposition 2.1.
This contradiction finishes the proof.
\end{proof}
\end{prop}

\begin{lem}
If $L$ is a line in $\P^4$ such that $L\cap \indet (F)$ is
nonempty and is not a single reduced point, then $L$ is contained
in $H$.
\end{lem}
\begin{proof}
Suppose $L$ intersects $\indet(F)$.  Then $L\cap \indet (F)$ is
zero dimensional, so its structure sheaf has finite dimension
$\lambda$ over the ground field $k$.  So the rational map $F|_L$
is given by $\O_{\P^1}(3-\lambda)$ after $F|_L$ is extended over
the points of indeterminacy.  The intersection $L\cap \indet (F)$
is a single reduced point exactly when $\lambda=1$.

If $\lambda =2$, then $F$ maps $L$ isomorphically onto a line.
This is impossible if $L$ is not contained in $H$ because $L$
would meet $X$ in a scheme of length 5 while $F(L)$ would meet $Y$
in a scheme of length 3.  If $\lambda=3$, then $F(L)$ is a point
in $Y$, so $L\subset H$.  This completes the proof because
$\lambda\leq 3$.
\end{proof}

\begin{prop}
$H$ is not $2K_1+2K_2$ for distinct hyperplanes $K_1$ and $K_2$.
\begin{proof}
Suppose $H=2K_1+2K_2$.  According to the formula of Corollary
3.10, there are points $p_i\in K_i$ such that $F|_{K_i}$ factors
through projection from $p_i$, and the tangent space to the
indeterminacy scheme $\indet(F|_{K_i})$ at $p_i$ is equal to the
tangent space of $K_i$ at $p_i$.

We claim $p_1\neq p_2$.  Indeed, if $p=p_1=p_2$, then
$T_{p}K_1\neq T_{p}K_2$ because the $K_i$ are distinct
hyperplanes. So $\indet(F)$ would have a four dimensional tangent
space at $p$.  But now every line in $\P^4$ through $p$ meets
$\indet(F)$ in a scheme that is nonreduced at $p$.  So Lemma 3.12
implies that every line in $\P^4$ through $p$ is contained in $H$.
This is impossible because $H$ is not all of $\P^4$, so $p_1\neq
p_2$.

Let $L$ be the line containing $p_1, p_2$. By Lemma 3.5, $H$ has
order at least 3 at both $p_i$, so $p_1, p_2\in K_1\cap K_2$.
Therefore, $L\subset K_1\cap K_2$, so that $L$ meets $\indet(F)$
in a scheme of length 6, which is impossible because $L\cap
\indet(F)$ can have length at most 3.
\end{proof}
\end{prop}

\begin{prop}
$H$ does not contain a hyperplane.
\begin{proof}
The only case left to rule out is $H=4K$.  So suppose $H=4K$.  We
will lift the polynomials that give $F|_K$ from $K$ to its second
infinitesimal neighborhood $2K$ and then derive a contradiction.

Using Corollary 3.10, we choose homogeneous coordinates
$x_0,\dots,x_4$ on $\P^4$ such that $K=V(x_4)$ and
 \begin{equation}
 F=(x_0^3+x_4q_0,\ x_1^3+x_4q_1,\ x_2^3+x_4q_2,\ x_0x_1x_2+x_4q_3,\
 x_4q_4)
 \end{equation}
for some homogeneous polynomials $q_i$ of degree 2 in the $x_i$.
Let $y_0,\dots,y_4$ be homogeneous coordinates on the target
$\P^4$.   Since $Y\cap Y^\prime$ is the image of $K=V(x_4)$ we see
$Y^\prime=V(y_4)$, and $Y\cap Y^\prime=V(y_4,\ y_0y_1y_2-y_3^3)$.
So $Y$ is given by the equation
 $$G:=y_0y_1y_2-y_3^3+y_4G_2$$
for some homogeneous $G_2$ of degree 2 in the $y_i$.

Since $2K\subset H$, $x_4^2$ divides the pullback
\begin{align*}
F^*G & \  = \ (x_0^3+x_4q_0)(x_1^3+x_4q_1)(x_2^3+x_4q_2)             \\
     & \ \ \ \ \ \ \ \ \    -(x_0x_1x_2+x_4q_3)^3+x_4q_4F^*G_2                              \\
    & \ = \
    x_4(q_0x_1^3x_2^3+q_1x_0^3x_2^3+q_2x_0^3x_1^3-3q_3x_0^2x_1^2x_2^2+q_4F^*G_2)  \\
    & \ \ \ \ \ \ \ \ \ +x_4^2 (\text{other terms})
\end{align*}
Therefore $x_4$ divides
\begin{equation}
q_0x_1^3x_2^3+q_1x_0^3x_2^3+q_2x_0^3x_1^3-3q_3x_0^2x_1^2x_2^2+q_4F^*G_2
\end{equation}
Using the notation $\partial_i:=\frac{\partial}{\partial y_i}$,
the partial derivatives of $G$ are
\begin{align*}
\partial_0G &= y_1y_2+y_4\partial_0G_2      \\
\partial_1G &= y_0y_2+y_4\partial_1G_2      \\
\partial_2G &= y_0y_1+y_4\partial_2G_2      \\
\partial_3G &= -3y_3^2+y_4\partial_3G_2      \\
\partial_4G &= G_2+y_4\partial_4G_2
\end{align*}
It follows that $V(y_0y_1,\ y_0y_2,\ y_1y_2,\ y_3,\ y_4,\ G_2)$ is
empty because it is contained in the singular locus of $Y$.  Hence
$V(G_2)$ does not contain the point $(1:0:0:0:0)$. In other words,
$y_0^2$ appears with nonzero coefficient in $G_2$. Therefore by
equation (3.6), $x_0^6$ has nonzero coefficient in $F^*G_2$.

Moreover, $q_4 F^*G_2$ is the only term in formula (3.7) in which
$x_0^6$ can occur because the $q_i$ have degree 2. Since (3.7) is
zero (mod $x_4$), the $x_0^6$ term in (3.7) disappears when
considered (mod $x_4$).  This can only happen if $x_4$ divides
$q_4$. Further consideration of (3.7) shows that there are scalars
$a_0, a_1, a_2$ such that the following equations hold (mod
$x_4$):
\begin{align*}
 q_0 &\equiv a_0x_0^2  \\
 q_1 &\equiv a_1x_1^2  \\
 q_2 &\equiv a_2x_2^2   \\
 q_3 &\equiv \frac{1}{3}\ (a_0x_1x_2+a_1x_0x_2+a_2x_0x_1)    \\
 q_4 &\equiv 0
\end{align*}
Therefore, equation (3.6) yields
\begin{align*}
F_0 &= x_0^3+a_0x_4x_0^2+x_4^2h_0    \\
F_1 &= x_1^3+a_1x_4x_1^2+x_4^2h_1    \\
F_2 &= x_2^3+a_2x_4x_2^2+x_4^2h_2    \\
F_3 &= x_0x_1x_2+\frac{1}{3}x_4(a_0x_1x_2+a_1x_0x_2+a_2x_0x_1)+x_4^2h_3  \\
F_4 &= x_4^2h_4
\end{align*}
for some homogeneous linear polynomials $h_i$ in $x_0,\dots,x_4$.

As $H=4K$ is contained in $F^{-1}(Y)$, $x_4^4$ divides
 $$ F^*G=F_0F_1F_2-F_3^3+F_4F^*G_2 .$$
Since $F^*G$ has degree 9, the monomial $x_0^6x_4^4$ can not
appear with nonzero coefficient in $F^*G$.  On the other hand, the
monomial $x_0^6$ does have nonzero coefficient in $F^*G_2$.  So
$x_0^6$ must appear in a monomial with nonzero coefficient in
 $$ F^*G - F_4F^*G_2 =F_0F_1F_2-F_3^3 . $$
However, by the formulas for $F_0,...,F_3$ given above, the
highest power of $x_0$ that can appear is $x_0^5$.  This
contradiction finishes the proof.
\end{proof}
\end{prop}

\begin{lem}
Every point $p\in\indet(F)$ is a point of order 3 on $H$.
\begin{proof}
Let $T_p$ denote the tangent space of the indeterminacy scheme
$\indet(F)$ at some closed point $p\in \indet(F)$. If $L$ is a
line in $\P^4$ tangent to $\indet(F)$ at $p$, then $L$ is
contained in $H$ by Lemma 3.12. Proposition 3.14 asserts that $H$
does not contain a hyperplane, so $T_p$ can have dimension at most
2.

By Lemma 3.5,  $H$ has order at least 3 at $p$.  So we only have
to rule out the case $\ord_pH=4$.

Suppose $p\in\indet(F)$ is a point of order 4 in $H$, so that $H$
is a cone over $p$.  Let $\pi:\widetilde \P^4\rightarrow \P^4$ be
the blowup of $\P^4$ at $\indet(F)$, and let $\widetilde
F:\widetilde\P^4\rightarrow Y$ be the resulting morphism that
extends $F$. Let $E_p=\pi^{-1} (p)$ be the preimage of the reduced
point $p$. Since the dimension of $T_p$ is at most 2, if $l_1$ and
$l_2$ are two general lines in $\P^4$ through $p$, then the scheme
$l_1\cup l_2$, which is contained in a 2-plane, meets $\indet(F)$
in the reduced point $p$. So the strict transform of $l_1\cup l_2$
is the blowup of $l_1\cup l_2$ at $p$, which is the
scheme-theoretic intersection $l_1\cap l_2$. Therefore the strict
transform of $l_1\cup l_2$ is the disjoint union of the strict
transforms $\tilde l_1$ and $\tilde l_2$ of $l_1$ and $l_2$. In
particular, a general line $l$ in $\P^4$ through $p$ determines a
point on $E_p$, and distinct general lines yield distinct points
in $E_p$. Hence $E_p$ has an irreducible component $E_p^\prime$
that is the birational image of the space $\P^3$ of lines in
$\P^4$ through $p$.

If $\tilde l$ is the strict transform of a general line $l$ in
$\P^4$ through $p$, then $\tilde l$ meets $\widetilde F^{-1}(Y)$
in a scheme of length 6.  Indeed, the morphism $\widetilde
F|_{\tilde l}$ is given by sections of the line bundle
$\O_{\P^1}(2)$ because the rational map $F|_l$ is given by
$\O_{\P^1}(3)$ and has indeterminacy scheme equal to a single
reduced point.  When the indeterminacy is resolved the resulting
line bundle is $\O_{\P^1}(2)$.  The Cartier divisor $\widetilde
F^{-1}(Y)$ is the pullback of the degree 3 divisor $Y$.  So
$\widetilde F^{-1}(Y)$ has intersection product 6 with the curve
$\tilde l$.

Moreover, $E_p^\prime$ maps onto $Y$, as follows.  If $\widetilde
H$ is the strict transform of $H$ in $\widetilde\P^4$, then
$\widetilde H$ intersects $E_p$ in a scheme of dimension 2.  In
particular, $E_p^\prime$ is not contained in $E_p\cap \widetilde
H$. For $l$ a general line in $\P^4$ through $p$, $l$ only meets
$H$ at $p$ because $H$ is a cone over $p$.  So the strict
transform $\tilde l$ of $l$ does not meet $\widetilde H$.

Since $\tilde l$ meets the strict transform $\widetilde X$ of $X$
in a scheme of length 5 and $\tilde{l}$ does not meet $\widetilde
H$, the other point of $\tilde l\cap \widetilde F^{-1}(Y)$ lies in
$E_p$ and therefore in $E_p^\prime$. So $\widetilde F$ maps
$E_p^\prime$ isomorphically onto $Y$ because $\widetilde F$ gives
a closed embedding of $E_p$ into $\P^4$ by Lemma 3.4.  This is a
contradiction because $E^\prime_p$ is rational and $Y$ is not.
Therefore every point $p\in \indet(F)$ is a point of order 3.
\end{proof}
\end{lem}

\begin{prop}
$H$ is not $2Q$ for an integral quadric $Q$.
\begin{proof}
By \cite[Proposition 2.2(2)]{Sheppard}, $Y$ is not the image of a
morphism from a hypersurface in $\P^4$ of degree 2, so $F|_Q$ is
not a morphism. Let $p\in \indet{F}$ be a closed point. If $H=2Q$,
then every closed point in $Q$ has order 2 or 4 in $H$.  By Lemma
3.5, $p$ is a point on $H$ of order at least 3.  So $p$ has order
4, which is impossible by Lemma 3.15.  We conclude $H\neq 2Q$.
\end{proof}
\end{prop}

\begin{prop}
$H$ is not $Q_1+Q_2$ for distinct integral quadrics $Q_1$ and
$Q_2$.
\begin{proof}
Suppose $H=Q_1+Q_2$.  By Lemma 3.8, the $Q_i$ map dominantly onto
$Y$. Therefore, by Lemma 3.3 applied to the blowups of $Q_1$ and
$Q_2$ at the indeterminacy schemes of $F|_{Q_1}$ and $F|_{Q_2}$,
for $L$ a general line in $Y$, $C:=F^{-1}(L)\cap H$ is reduced at
the generic point of each of its irreducible components.

By Lemma \cite[Proposition 2.2(2)]{Sheppard}, $F|_{Q_1}$ can not
be a morphism because $\deg Y>2$. Pick $p\in\indet(F|_{ Q_1})$. By
Lemma 3.15, $p$ is a point of order 3 in $H$, so $Q_2$ contains
$p$. I claim that $C$ has an irreducible component contained in
each of the $Q_i$ through $p$. Indeed, if $\pi_i:\widetilde
Q_i\rightarrow Q_i$ is the blowup at $\indet(F|_{Q_i})$, and $q_i:
\widetilde Q_i\rightarrow Y$ is the resulting morphism, then a
general line $L$ in $Y$ will meet $q_i(\pi_i^{-1}(p))$. So the
preimage of $L$ in $Q_i$ passes through $p$.

Therefore, by Corollary 3.7, the connected component of $C$
containing $p$ is the union of two lines through $p$, one line
contained in $Q_1$ and the other in $Q_2$. But now each of the
$Q_i$ contain a two parameter family of lines through $p$ because
$Y$ is covered by a two parameter family of lines. So both the
$Q_i$ are cones over $p$, so that $p$ is a point of order 4 in
$H$.  This contradicts Lemma 3.15.
\end{proof}
\end{prop}

The following Lemma will be needed to analyze the preimage in $H$
of a general line in $Y$ in the case where $H$ is integral.

\begin{lem}
Let $\pi:\mathbb F\rightarrow B$ be a projective morphism with $B$
integral. Suppose $\sigma:B\rightarrow \mathbb F$ is a section of
$\pi$, and for each $b\in B$ the connected component $\mathbb F_b$
of $\pi^{-1}(b)$ that contains $\sigma(b)$ is irreducible. Then
$\bigcup_b \mathbb F_b$ is an irreducible component of $\mathbb
F$.
\begin{proof}
Let $\mathbb
F\stackrel{\pi^\prime}{\rightarrow}B^\prime\stackrel{g}{\rightarrow}B$
be the Stein factorization of $\pi$.  So $\pi^\prime$ has
connected fibers, and $g$ is finite.  Then $\pi^\prime\circ\sigma$
is a section of $g$, and
$(\pi^\prime)^{-1}(\pi^\prime\circ\sigma)(b)=\mathbb F_b$. Notice
that $(\pi^\prime\circ\sigma)(B)$ is an irreducible component of
$B^\prime$ because they have the same dimension.  Since
$(\pi^\prime)^{-1}(\pi^\prime\circ \sigma(B))$ has irreducible
fibers $\mathbb F_b$ and $B$ is irreducible,
$(\pi^\prime)^{-1}(\pi^\prime\circ \sigma(B))=\bigcup_b \mathbb
F_b$ is irreducible.
\end{proof}
\end{lem}

The only case left to rule out is when $H$ is integral.

\begin{thm}
$H$ is not integral.
\begin{proof}
Suppose $H$ is integral.  By Lemma 3.8, $F|_H$ is dominant because
$H$ has multiplicity 1 in $F^{-1}(Y)$. However, $F|_H$ can not be
a morphism, because if it were then Table 1 in the Appendix would
imply that the polynomial degree $m$ of $F$ would be at most 2,
not 3.  Choose $p\in\indet(F)$ a reduced point.

Let $L$ be a general line in $Y$, cut out by the linear forms
$\xi_1, \xi_2, \xi_3$ in $\P^4$.  Define
 \begin{align*}
  F^{-1}(L) &:=V(F^*\xi_1,\ F^*\xi_2,\ F^*\xi_3) \\
  C(L)&:=F^{-1}(L)\cap H \\
  D(L)&:=F^{-1}(L)\cap X
 \end{align*}
By Lemma  3.3, $C=C(L)$ and $D=D(L)$ are reduced at the generic
point of each of their irreducible components, and

Since the linear forms $F^*\xi_i$ vanish on the indeterminacy
scheme $\indet(F):=V(F_0,\dots, F_4)$, $\indet(F)$ is contained in
$F^{-1}(L)=C\cup D$.  Therefore, $\indet(F)$ is contained in $C$
because $\indet(F)$ does not intersect $D$, as $D$ is contained in
$X$.  By Corollary 3.7, the connected component of $C$ that
contains $p$ is either a smooth quadric curve or the union of two
lines meeting at some point.  According to Lemma 3.15, $p$ is
triple point on $H$, so that $H$ is not a cone over $p$.  By
Corollary 3.7, the connected component of $C$ that contains $p$ is
a plane conic.  Since $H$ is not a cone over $p$, $Y$ has a two
dimensional family of lines, and $L$ is a general line on $Y$, the
connected component of $C$ that contains $p$ is a smooth plane
conic.

Let $B$ be an open subscheme of the space of lines in $Y$ such
that for every $L\in B$, every component of $C(L)$ that meets
$\indet(F)$ is a smooth plane conic.  For $L\in B$, $F^{-1}(L)$
has degree 27, and $D(L)$ has degree 15.  So $C(L)$ has degree 12
and is the disjoint union of 6 reduced plane conics.

Choose a general $L_0\in B$, and let $C_1(L_0),\dots, C_6(L_0)$ be
the connected components of $C(L_0)$.  Let $\zeta_i :=
C_i(L_0)\cap \indet(F)$, and let $\lambda(\zeta_i)$ denote the
length of the zero dimensional scheme $\zeta_i$.  The polynomial
degree of $F$ is $m=3$, $C_i(L_0)$ has degree 2, and $L_0$ has
degree 1.  So the restriction $C_i(L_0)\dashrightarrow L_0$ of $F$
to $C_i(L_0)$ has degree $6-\lambda(\zeta_i)$.  Therefore $\deg
F|_H$ is the sum of the $\deg F|_{C_i(L_0)}$:
\begin{equation}
  \deg F|_H = \sum_{i=1}^6 6-\lambda(\zeta_i).
\end{equation}
Now number the $C_i(L_0)$ so that $\lambda(\zeta_1)\geq
\dots\geq\lambda(\zeta_6)$, and let $p$ be a closed point of
$\zeta_1$.  Let $\mathbb F\subset B\times H$ be the total space of
the family $\pi:\mathbb F\rightarrow B$ whose fiber over $L\in B$
is $C(L)$.  The closed subscheme $B\times p\subset \mathbb F$ is a
section of $\pi$.  So
 $$ \mathbb F_1 := \bigcup_{L\in B}C_1(L) $$
is an irreducible component of $\mathbb F$ by Lemma 3.18.

Notice that there is a morphism $\mathbb F\rightarrow B_H$, where
$B_H$ is the space of quadric plane curves in $H$ that contain
$p$, given by sending a point $ (L, x)\in \mathbb F\subset B\times
H$ to the irreducible component of $C(L)$ that contains $x$. The
fibers of $\mathbb F\rightarrow B_H$ are one dimensional and $\dim
\mathbb F =3$.  So the image in $B_H$ of the intersections of the
various irreducible components of $\mathbb F$ has dimension at
most 1.  Hence there is at most a one dimensional space of lines
in $Y$ such that $C(L)$ has a component corresponding to a point
in $B_H$ whose fiber in $\mathbb F$ lies in more than one
irreducible component. So by generality of $L_0\in B$, $\mathbb
F_1$ is the only irreducible component of $\mathbb F$ that has
nonempty fiber over the point $C_1(L_0)\in B_H$. Therefore, the
irreducible component $\mathbb F_1$ of $\mathbb F$ did not depend
on the choice of $p\in\zeta_1$ because for any $p\in\zeta_1$ the
connected component of $C(L_0)$ that contains $p$ is $C_1(L_0)$.

Since $\mathbb F_1$ did not depend on the choice of $p\in
\zeta_1$, for every line $L\in B$, the component
$C_1(L):=\pi^{-1}([L])\cap\mathbb F_1$ of $C(L)$ that is contained
in $\mathbb F_1$ has the property that $C_1(L)\cap
\indet(F)=\zeta_1$.  This is because for every $q\in\zeta_1$,
$B\times q$ is contained in $\mathbb F_1\subset B\times H$.

Consider the composite morphism $\phi:\mathbb F_1 \rightarrow H
\stackrel {F|_H}{\dashrightarrow} Y$. Let $y\in Y$ be a general
point. There are six lines $L_1, \dots, L_6$ on $Y$ through $y$,
and $\phi^{-1}(L_i)=C_1(L_i)$.  By counting preimage points of $y$
we see that $\deg \phi$ is the sum of the degrees of the
$C_1(L_i)\dashrightarrow L_i$:
\begin{equation}
\deg \phi = \sum_{i=6}^6 6-\lambda(\zeta_1).
\end{equation}
Since $\deg \phi\geq\deg F|_H$, equations (3.8) and (3.9) show
that all the $\lambda(\zeta_i)$ are equal by maximality of
$\lambda(\zeta_1)$. Therefore $\mathbb F_1\rightarrow H$ is a
birational morphism.

Let $p\in\indet(F)$.  By Lemma 3.15, $p$ is a triple point on $H$,
and so $H$ is rational.  Therefore $\mathbb F_1$ is rational and
dominates the surface $B$. However, $B$ does not contain a
rational curve by Proposition 2.1. This contradiction shows that
$H$ can not be integral.
\end{proof}
\end{thm}

% ------------------------------------------------------------------------
%GATHER{Xbib.bib}   % For Gather Purpose Only
%GATHER{Paper.bbl}  % For Gather Purpose Only
\bibliographystyle{amsplain}
\bibliography{xbib}
\end{document}